\documentclass[10pt]{amsart}
\usepackage{amssymb}
\usepackage{amsmath,amssymb, amsthm}

\newcommand{\epss}{\varepsilon}

\newcommand{\SL}{\mathrm{SL}(2,\mathbb{R})}

\def\tA{\mathcal{A}}
\def\tB{\mathcal {B}}
\def\tC{\mathcal {C}}

\def\a{\alpha}
\def\w{\omega}

\newcommand{\bnu}{{\bar\nu}}

\let\newpf\proof \let\proof\relax 
\newenvironment{pf}{\newpf[\proofname]}{\qed\endtrivlist}

\newcommand{\norm}[1]{\bigl\| #1 \bigr\|}
\renewcommand{\norm}[1]{|\!|\!| #1 |\!|\!|}

\newcommand{\bQ}{\overline{Q}}
\newcommand{\ba}{\overline{A}}

\newcommand{\btau}{\overline{\tau}}

\newcommand{\cO}{\mathcal{O}}

\newcommand{\CD}{\rm CD}

\def\be{\begin{equation}}
\def\ee{\end{equation}}

\def\ba{{\begin{align}}}
\def\ea{{\end{align}}}

\def\bm{\begin{matrix}}
\def\em{\end{matrix}}

\renewcommand{\sl}{{\mathrm{sl}}}

\def\a{{\alpha}}

\def\tA{\tilde {\mathcal {A}}}
\def\tB{\tilde {\mathcal {B}}}
\def\tC{\tilde {\mathcal {C}}}

\def\SL{{\mathrm{SL}}}
\def\PSL{{\mathrm{PSL}}}
\def\SO{{\mathrm{SO}}}

\def\0{{\mathbf 0}}

\def\cal{\mathcal}

\newtheorem{thm}{Theorem}[section]
\newtheorem{cor}[thm]{Corollary}

\newtheorem{lem}[thm]{Lemma}
\newtheorem{lemma}[thm]{Lemma}
\newtheorem{claim}[thm]{Claim}
\newtheorem{prop}[thm]{Proposition}

\theoremstyle{remark}
\newtheorem{rem}{Remark}[section]

\newtheorem{problem}{Problem}

\numberwithin{equation}{section}

\def \bn {\hfill \\ \smallskip\noindent}

\theoremstyle{definition}
\newtheorem{defn}{Definition}[section]

\def\proof{\bn {\bf Proof.} }

\def\epss{\varepsilon}

\def\ssm{\smallsetminus}

\renewcommand{\setminus}{\ssm}

\newcommand{\id}{\operatorname{id}}

\newcommand{\eps}{{\epsilon}}

\newcommand{\QQ}{{\cal Q}}

\newcommand{\C}{{\mathbb C}}

\newcommand{\M}{{\mathbb M}}
\newcommand{\N}{{\mathbb N}}
\newcommand{\Q}{{\mathbb Q}}
\newcommand{\R}{{\mathbb R}}
\newcommand{\T}{{\mathbb T}}

\newcommand{\Z}{{\mathbb Z}}

\def\B0{{\bold{0}}}

\catcode`\@=12

\def\Empty{}
\newcommand\oplabel[1]{
  \def\OpArg{#1} \ifx \OpArg\Empty {} \else
  	\label{#1}
  \fi}
		
\newcommand{\comm}[1]{}
\newcommand{\comment}[1]{}

\begin{document}

\title[Liouvillean KAM]{A KAM scheme for $\SL(2,\R)$
cocycles with Liouvillean frequencies}

\author {Artur Avila}
\address{
CNRS UMR 7586, Institut de Math\'ematiques de Jussieu\\
175, rue du Chevaleret\\
75013--Paris, France
}
\email{artur@math.sunysb.edu}

\author{Bassam Fayad}
\address{
CNRS UMR 7539, LAGA\\
Universit\'e Paris 13\\
93430--Villetaneuse, France
}
\email{fayadb@math.univ-paris13.fr}

\author{Rapha\"el Krikorian}

\address{
CNRS UMR 7599, Laboratoire de Probabilit\'es et Mod\`eles al\'eatoires\\
Universit\'e Pierre et Marie Curie--Boite courrier 188\\
75252--Paris Cedex 05, France
}
\email{raphael.krikorian@upmc.fr}

\date{\today}

\begin{abstract}

We develop a new KAM scheme that applies to $\SL(2,\R)$ cocycles with one
frequency, irrespective of any Diophantine condition on the base dynamics.
It gives a generalization of Dinaburg-Sinai's Theorem to arbitrary
frequencies: under a closeness to constant assumption,
the non-Abelian part of the classical reducibility problem can
always be solved for a positive measure set of parameters.

\end{abstract}

\setcounter{tocdepth}{1}

\maketitle

\section{Introduction}

In this paper we are concerned with analytic
quasiperiodic $\SL(2,\R)$ cocycles in
one frequency.  Those are linear skew-products
\begin{align}
(\alpha,A):&\T \times \R^2 \to \T \times \R^2\\
\nonumber
&(x,w) \mapsto (x+\alpha,A(x) \cdot w),
\end{align}
where $\alpha \in \R$ and $A:\T \to \SL(2,\R)$ is analytic ($\T=\R/\Z$).
The main source of examples are given by Schr\"odinger cocycles, where
\be
A(x)=S_{v,E}(x)=\left (\bm E-v(x) & -1 \\ 1 & 0 \em \right ),
\ee
which are related to one-dimensional quasiperiodic Schr\"odinger operators
\be
(Hu)_n=u_{n+1}+u_{n-1}+v(\theta+n \alpha) u_n.
\ee

We are interested in the case where $A$ is close to a constant.  In this
case, the classical question is whether $(\alpha,A)$ is reducible.  This
means that $(\alpha,A)$ is conjugate to a
constant, that is, there exists $B:\T \to \PSL(2,\R)$ analytic such that
$B(x+\alpha)A(x)B(x)^{-1}$ is a constant.

Reducibility could be thought as breaking up into two different problems. 
One could first try to conjugate the cocycle into some Abelian subgroup of
$\SL(2,\R)$, and then conjugate to a constant inside the subgroup.  It is
easy to see that the second part involves the solution of a cohomological
equation.  Thus small divisor obstructions related to $\alpha$ must
necessarily be present in the problem of reducibility.
However, if one takes the point of view that finding a solution of the cohomological
equation is an understood problem,\footnote {We must point out that the
solution of the cohomological equation is an interesting problem
in some concrete situations, see \cite {AJ1}.}
we should shift our focus to understanding the first problem.

We focus on the case where the Abelian subgroup is $\SO(2,\R)$.  Let us say
that $(\alpha,A)$ is conjugate to a cocycle of rotations if there exists
$B:\T \to \SL(2,\R)$ analytic such that $B(x+\alpha)A(x)B(x)^{-1} \in
\SO(2,\R)$.  We obtain the somewhat surprising result that conjugacy to
a cocycle of rotations is frequent (under a closeness to constant
assumption) irrespective of any condition on $\alpha$.  A particular
case that ilustrates what we mean by frequent is the following:

\begin{thm} \label{dsliouville}

Let $v:\T \to \R$ be analytic and close to a constant.
For every $\alpha \in \R$, there exists a positive measure set of $E \in
\R$ such that $(\alpha,S_{v,E})$ is conjugate to a cocycle of
rotations.\footnote {Indeed the Lebesgue measure of the set $X(v,\alpha)$
of energies such that $(\alpha,S_{v,E})$ is conjugate to a cocycle of
rotations converges (uniformly on $\alpha$) to $4$ as $v$ converges to a
constant.  We recall that for $v$ constant, $X(v,\alpha)$
is an open interval of length $4$, and for $v$ non-constant $X(v,\alpha)$
(and indeed the larger set of energies where the Lyapunov exponent
vanishes) has Lebesgue measure strictly less than $4$.}

\end{thm}

\begin{rem}

The existence of a positive measure set of energies for which
$(\alpha,S_{v,E})$ is conjugate to a cocycle of rotations
implies the existence of some absolutely continuous spectrum for the
corresponding Schr\"odinger operator (for stronger results, see
\cite {LS}).  In this context and under
the additional condition that $v$ should be {\it even}, the
existence of some absolutely continuous spectrum was first
obtained by Yoram Last (unpublished), by different methods.

\end{rem}

Another consequence of our
new approach to the questions of reducibility is the
following generalization of the global almost sure dichotomy 
between non uniform hyperbolicity and reducibility of Schr\"odinger
cocycles obtained in \cite{AK1} for recurrent
Diophantine frequencies in the base (and extended in \cite{AJ2,FK}
to all Diophantine frequencies).

\begin{thm} \label{global} Let $\a \in \R \setminus \Q$ and $v:\T \to \R$ be analytic, then for almost every  
energy $E \in
\R$ the cocycle $(\alpha,S_{v,E})$ is either non-uniformly hyperbolic or (analytically) conjugate to a cocycle of rotations.
\end{thm}

\begin{rem}

Theorem \ref {global} has the following interesting consequence:
in the essential support of the absolutely continuous spectrum,
the generalized eigenfunctions are almost surely bounded.  This property is
in fact expected to hold for general ergodic Schr\"odinger operators (the
so-called Schr\"odinger Conjecture \cite {MMG}), but this is the first time
it is verified in a Liouvillean context.

\end{rem}

The proof of Theorem \ref{global} starts from Kotani's theory that
essentially asserts that there is an almost sure dichotomy between
non-uniform hyperbolicity and $L^2$ rotations-reducibility
(conjugacy to a cocycle with values in $\SO(2,\R)$ of the
cocycles  $(\alpha,S_{v,E})$. Then it uses the convergence to constants of the  renormalizations  of  cocycles that are $L^2$-conjugated to rotations, obtained in \cite{K, AK1,AK2}.
Finally, it is concluded using a precise version of the local result given in Theorem \ref{dsliouville}, namely Theorem \ref{theorem.cd} below.

The proof of Theorem \ref{dsliouville}  will be based on a new KAM scheme for $\SL(2,\R)$ cocycles
that is able to bypass the small divisors related to $\alpha$.  An important
ingredient is the ``cheap trick'' recently developed in \cite {FK}, though
our implementation is rather different (it does not involve renormalization,
and it is applied without any Liouvillean assumption on $\alpha$).

Before giving a more general statement, let us first recall some
of the history of KAM methods applied to reducibility problems.

The first result of reducibility was due to Dinaburg-Sinai \cite {DS}.
It establishes that $(\alpha,S_{v,E})$ is reducible for a positive measure
set $E$ whenever $\alpha$ is Diophantine and $v$ is sufficiently close to a
constant (depending on $\alpha$).

\comm{
In certain parametrized families of cocycles $(\alpha,A_E)$,
with $\alpha$ Diophantine and the
$A_E$ sufficiently close to constant (depending on the Diophantine
condition), there exists a positive measure set of parameters $E$ for which
$(\alpha,A_E)$ is reducible.
}

A more precise formulation (due to
Herman \cite {H}), gives an explicit condition on the topological dynamics
of $(\alpha,A)$ which guarantees the reducibility (under closeness to
constant assumptions).  It is formulated in terms of the fibered rotation
number $\rho=\rho(\alpha,A) \in \T$,
a topological invariant which can be defined for all
cocycles homotopic to the identity (this includes all Schr\"odinger cocycles
and all cocycles which are close to constant).  To define this invariant,
first notice that
$(\alpha,A)$ naturally acts on the torus $\T \times \T$ if one
identifies the second coordinate with the set of arguments
of non-zero vectors of $\R^2$.  This torus map turns out to have a well
defined rotation vector (if $\alpha \in \R \setminus \Q$, otherwise one
averages over the first coordinate).
The first coordinate of the rotation vector is necessarily $\alpha$ and
the second is, by definition the fibered rotation number $\rho$.  As a
Corollary of the classical KAM Theorem on $\R^2/\Z^2$, Herman concludes
that if $(\alpha,\rho)$ is Diophantine and if $A$ is sufficiently close to
a constant (depending on $\alpha$ and $\rho$), then $(\alpha,A)$ is
reducible.

\comm{
This Theorem can be proved by a conventional KAM method (actually it can be
seen as a corollary of the classical KAM theorem on $\R^2/\Z^2$).
}

Herman's Theorem suggests that, in the theory of reducibility, there are
obstructions (small divisors) related not only to $\alpha$, but also
on the joint properties of the pair $(\alpha,\rho)$.  That this is actually
the case is not obvious.  It follows from a result of Eliasson
that for every $\alpha$ there exists a generic set $B_\alpha$ such that for
every $\rho \in B_\alpha$ there exists $A$ arbitrarily
close to constant with $\rho(\alpha,A)=\rho$, and such that $(\alpha,A)$ is
not reducible.\footnote {\label {js}
More explicit examples can be given from the theory of the almost Mathieu
operator.  Let $v(x)=2 \cos 2 \pi x$,
let $\alpha \in \R$ be arbitrary, and let
$\rho \in [0,1/2]$ be such that
for every $k \in \Z$ $\|2 \rho-k \alpha\|_\T>0$, while
$\liminf_{|k| \to \infty} \frac {1} {|k|}
\ln \|2 \rho-k\alpha\|_\T=-\infty$.
It is well known that
for every $\lambda \neq 0$, there exists a unique
$E=E(\lambda,\alpha,\rho) \in \R$
such that $\rho(\alpha,S_{\lambda v,E})=\rho$.  On the other hand, it can be
shown that $(\alpha,S_{\lambda v,E})$ is not reducible: indeed reducibility
implies
the existence of an eigenfunction for the ``dual'' operator
$(H u)_n=u_{n+1}+u_{n-1}+\lambda^{-1} v (\rho+n \alpha) u_n$
with energy $\lambda^{-1} E$ (this is an instance of Aubry duality),
but $H$ has no point spectrum by \cite {JS}.}

\def\cQ{{\mathcal{Q}}}
\def\cP{{\mathcal{P}}}

We can now state the precise version of our Main Theorem.  Let
$\Delta_h=\{x \in \C/\Z,\, |\Im x|<h\}$.  For a bounded
holomorphic (possibly matrix
valued) function $\phi$ on $\Delta_h$, we let $\|\phi\|_h=\sup_{x \in
\Delta_h} \|\phi(x)\|$. We denote by $C^\w_{h}(\T,*)$ the set of all these
$*$-valued functions ($*$ will usually denote $\R$, $\SL(2,\R)$, or
the space of $2$-by-$2$ matrices $\M(2,\R)$).
In the theorem $\a$ is only supposed to be irrational and $q_n$ denotes the sequence 
of denominators of its best rational approximations (see section \ref{sec:2.3}).
For $\tau>0$, $0<\nu<1/2$ and $\epss>0$, let $\cQ_\alpha(\tau,\nu,\epss)
\subset \T$ be the set of all $\rho$ such that for every $i$ we have
\be
\|2 q_i \rho\|_{\T}>\epss \max \{q^{-\nu}_{i+1},q^{-\tau}_i\}.
\ee

\begin{thm}
\label{theorem.cd}
For every $\tau>0$, $0<\nu<1/2$, $\epss>0$, $h_*>0$, there exists
$\epsilon=\epsilon(\tau,\nu,\epss)>0$ with the following property.
Let $h>h_*$ and let
$A \in C^\omega_h(\T,\SL(2,\R))$ be real-symmetric such that
$\|A-R\|_{h}<\epsilon$ for some rotation matrix $R$, and
$\rho=\rho(\alpha,A) \in \cQ_\alpha(\tau,\nu,\epss)$.
Then there exists real-symmetric $B:\Delta_{h-h_*} \to {\rm SL}(2,\C)$
and $\phi:\Delta_{h-h_*} \to \C$ such that
\begin{enumerate}
\item $\|B-\id\|_{h-h_*}<\epss$ and
$\|\phi-\rho\|_{h-h_*}<\epss$,
\item $B(x+\alpha)A(x)B(x)^{-1}=R_{\phi(x)}$.
\end{enumerate}

\end{thm}

\noindent {\it Proof of Theorem \ref {dsliouville}.}
We may assume that $v$ has average $0$.
Recall (see, e.g., \cite {H}) that the function
$\rho:E \mapsto \rho(\alpha,S_{E,v})$ is continuous,
non-increasing and onto $[0,1/2]$.
Fix $\tau>0$, $0<\nu<1/2$ and
$\epss>0$ so that $X=\QQ_\alpha(\tau,\nu,\epss) \cap [0,1/2]$ has positive
measure for every $\alpha$, and let $Y=\rho^{-1}(X)$.
If $v$ is sufficiently close to constant
then Theorem \ref {theorem.cd} implies that
$(\alpha,S_{E,v})$ is analytically rotations reducible for every $E \in Y$.
Moreover, at any $E$ such that
$(\alpha,S_{E,v})$ is analytically rotations-reducible, it is easy to see
that the fibered rotation number has a bounded derivative.  Since $X$ has
positive Lebesgue measure and $\rho(Y)=X$, it follows that $Y$ has positive
Lebesgue measure as well.\footnote {In order to get that $Y$ has measure
close to $4$, one should use that for $E \in Y$,
Theorem \ref {theorem.cd} actually
guarantees the existence of a conjugacy to rotations
$B \in C^\omega(\T,\SL(2,\R))$ which is close to a constant matrix
conjugating $S_{0,E}$ to a rotation.  The estimate can then be concluded in
various ways, for instance as a consequence of the formula:
$\frac {d\rho} {dE}=-\frac {1} {8 \pi} \int_{\T} \|B\|^2_{\mathrm {HS}}$
(whose verification is straightforward).
}

\noindent {\it Proof of Theorem \ref{global}.}  We just sketch the argument,
which is a standard ``global to local'' reduction through convergence of
renormalization (previously established in \cite {AK2}).

Fix $\alpha \in (0,1) \setminus \Q$, and let
$\alpha_n \in (0,1) \setminus \Q$, $n \geq 0$, be the
$n$-th iterate of $\alpha$ by the Gauss map $x \mapsto \{x^{-1}\}$
(see section \ref{sec:2.3}).  Let $\beta_n=\prod_{k=0}^n \alpha_k$.

\begin{lemma} \label {ren}

Let $\cP \subset [0,1)$ be the set of all $\rho$
such that there exist $\tau>0$, $0<\nu<1/2$ and $\epss>0$ such that
$\frac {\rho} {\beta_{n-1}} \in \cQ_{\alpha_n}(\tau,\nu,\epss)$ for
infinitely many $n$.  Let $A \in
C^\omega(\T,\SL(2,\R))$ be homotopic to a constant and such that
$(\alpha,A)$ is $L^2$-conjugated to an $\mathrm{SO}(2,\R)$-valued cocycle
and the fibered rotation number of $(\alpha,A)$
belongs to $\cP$.
Then $(\alpha,A)$ is analytically rotations-reducible.

\end{lemma}

\begin{pf}

Renormalization of cocycles (see \cite {K}, \cite {AK1})
associates to each analytic conjugacy class $[(\alpha,A)]$ of cocycles in
$(\R \setminus \Q) \times C^\omega(\T,\SL(2,\R))$, an analytic conjugacy
class $[(\alpha_n,A_n)]$, such that
if $(\a_n,A_n)$ is rotations-reducible (or reducible) then
$(\a,A)$ is rotations-reducible (or reducible) in the same class of
regularity.

If $A \in C^\omega(\T,\SL(2,\R))$ is homotopic to a constant and
$(\alpha,A)$ is
$L^2$-conjugated to an $\SO(2,\R)$-valued cocycle, then it follows from
convergence of renormalization \cite {AK2} that there are {\it
renormalization representatives} $(\alpha_n,A_n)$ such that
\begin{itemize}
\item the fibered rotation number of $(\alpha_n,A_n)$ is
$\rho_n=\frac {\rho} {\beta_{n-1}}$, where $\rho \in [0,1)$ is the fibered
rotation number of $(\alpha,A)$ and $\beta_{n-1}=\prod_{j=0}^{n-1}\a_j$,
\item there exists $h_0>0$ such that each $A_n$ admits a
bounded holomorphic extension to $\Delta_{h_0}$ and
$\|A_n-R_{\rho_n}\|_{h_0} \to 0$.
\end{itemize}

By definition of $\cP$, if the fibered rotation number of $(\alpha,A)$
belongs to $\cP$ then we can find $\tau>0$, $0<\nu<1/2$, $\epss>0$ and
arbitrarily large $n \geq 0$ such that
$\rho_n \in \cQ_{\alpha_n}(\tau,\nu,\epss')$.  Let $\epsilon=
\epsilon(\tau,\nu,\epss)$ be as in Theorem \ref {theorem.cd}, and choose such an
$n$ so large that
$\|A_n-R_{\rho_n}\|_h<\epsilon$.  By Theorem \ref {theorem.cd}, $(\alpha_n,A_n)$,
and hence $(\alpha,A)$,
must be analytically rotations-reducible.
\end{pf}

On the other hand, Kotani Theory (interpreted in the language of cocycles,
see \cite {AK1}, Section 2), yields:

\begin{lemma} \label {e}

Let $\cP \subset [0,1)$ be any full measure subset.  For every
$v \in C^\omega(\T,\R)$,
for almost every $E \in \R$, $A=S_{v,E}$ satisfies:
\begin{itemize}
\item either $(\alpha,A)$ has a positive Lyapunov exponent, or
\item $(\alpha,A)$ is
$L^2$-conjugated to an $\mathrm{SO}(2,\R)$-valued cocycle and the fibered
rotation number of $(\alpha,A)$ belongs to $\cP$.
\end{itemize}

\end{lemma}

The result follows from a combination of Lemmas \ref {ren} and \ref {e},
since the set $\cP$ specified in Lemma \ref {ren} has full measure
in $[0,1)$ by a simple application of Borel-Cantelli.
\qed

\subsection{Further remarks}

\subsubsection{Eliasson's Theory}

Eliasson developed a non-standard KAM scheme which
enabled him to prove a much stronger version of Dinaburg-Sinai's Theorem. 
He showed that for every Diophantine $\alpha$, there exists a full Lebesgue
measure set $\Lambda(\alpha)$ (explicitly given in terms of a Diophantine
condition)
\comm{
\be
\Lambda(\alpha)=\{\rho \in \T,\,
\limsup_{|k| \to \infty} \frac {1} {\ln |k|}
\ln \|2 \rho-k \alpha\|_{\T}>-\infty\}
\ee
}
such that if $A$ is sufficiently close
to constant (depending on $\alpha$) and if $\rho(\alpha,A) \in
\Lambda(\alpha)$ then $(\alpha,A)$ is reducible.
The natural question raised by our work is whether it is possible to
combine the strength of Eliasson's Theorem and that of our main result:

\begin{problem} \label {fullmeasure}

Let $\alpha \in \R \setminus \Q$\footnote {The result is
obviously false for $\alpha \in \Q$.}.  Is there a full
Lebesgue measure subset $\Lambda(\alpha)$ (explicitly given in terms of some
Diophantine condition) such that for every $A$ sufficiently close to
constant such that $\rho(\alpha,A) \in \Lambda(\alpha)$,
$(\alpha,A)$ is conjugate to a cocycle of rotations?

\end{problem}

Such a result would be
necessarily very subtly dependent on the analytic regularity of the cocycle
since there are recent counterexamples in Gevrey classes \cite {AKgevrey}.
For this reason, it is unlikely that a positive solution to the problem
could ever be achieved by a KAM method (since KAM methods
tend to generalize to weaker regularity).
On the other hand, \cite {AJ2} showed that the closeness
quantifier in Eliasson's Theorem was $\alpha$-independent, as long as
$\alpha$ was Diophantine 
which leads us to believe that it is possible to obtain a result for
all $\alpha$ in the analytic category.  Unfortunately, the method of \cite
{AJ2} mixes up at its core the Abelian and non-Abelian
subproblems of reducibility, so it seems unsuitable to prove results of the
kind obtained here.  A solution to Problem \ref {fullmeasure} has been obtained by the
first author \cite {A}: the proof provides (without restrictions on the fibered
rotation number) a sequence of conjugacies which put the
cocycle arbitrarily close to constants, so that (under a full measure
condition on the fibered rotation number) the results of this paper
eventually can be applied.

\subsubsection{Parabolic behavior}

Besides the elliptic subgroup $\SO(2,\R)$,
there are essentially two other Abelian subgroups of $\SL(2,\R)$.
Those are the hyperbolic subgroup
(diagonal matrices) and the parabolic subgroup
(stabilizer of a non-zero vector).
It has been understood for quite some time
that a relatively simple and open condition (uniform hyperbolicity of the cocycle)
ensures that the cocycle is conjugate 
 to a diagonal cocycle (this does not
involve any KAM scheme).  This condition is satisfied
for an open set of cocycles.  The parabolic case is, on the other hand, in
the frontier of elliptic and hyperbolic behavior and so it has ``positive
codimension'' (for instance,
it happens at most for a countable set of energies for
Schr\"odinger cocycles).  If $\alpha$ is Diophantine, Eliasson showed that
if $A$ is close to constant (depending on $\alpha$) then
$(\alpha,A)$ is conjugate  
degenerate) if and only if $(\alpha,A)$ is not uniformly hyperbolic and
$\|2 \rho(\alpha,A)-k \alpha\|_{\T}=0$ for some $k \in \Z$.

Unfortunately, no result of the kind obtained in this paper can
be proved in the parabolic case.\footnote {
For an explicit counterexample, consider $v(x)=2 \cos 2 \pi x$, and let $\alpha$ be
such that $\liminf_{q \to \infty} \frac {1} {q} \ln \|q
\alpha\|_{\T}=-\infty$.  Then for every $\lambda \neq 0$ and every $E \in
\R$, the cocycle $(\alpha,S_{\lambda v,E})$ is not conjugate
to a parabolic subgroup of $\SL(2,\R)$.  Indeed, as in footnote \ref {js},
such a conjugacy implies the existence of an eigenvalue for a dual
operator (by Aubry duality), which is known to have no point spectrum
(by Gordon's Lemma, \cite {AS}).}
In fact even the less ambitious goal of
analytically conjugating a cocycle to a triangular form can not be
obtained without some arithmetic restriction on $\alpha$: consider, for
instance, a cocycle of the form $(\alpha,A)$ with $A=R_{\lambda \phi}$, where
$\phi \in C^\omega(\T,\R)$ has average $0$ but is not a coboundary, i.e.,
the equation $\phi(x)=\psi(x+\alpha)-\psi(x)$ admits no solution $\psi \in
C^\omega(\T,\R)$ (such a $\phi$ exists as long as $\alpha$ can be
exponentially well approximated by rational numbers).  Then
$\rho(\alpha,A)=0$, but $(\alpha,A)$ can not be analytically conjugated to a
cocycle in triangular form (from such a conjugacy it is easy to construct a
solution of the cohomological equation).

\comm{
to a triangular
subgroup of $\SL(2,\R)$ can .
but it is unclear whether a result of
conjugacy to a triangular subgroup can be proved instead.
The following gives a possible formulation of such a result:

\begin{problem}

Let $\alpha \in \R \setminus \Q$ (this result is obviously false for $\alpha
\in \Q$), and let $\rho$ be such that $\|2 \rho-k
\alpha\|_{\T}=0$ for some $k \in \Z$.  Is it true that if $A$ is close to
constant and $\rho(\alpha,A)=\rho$ then there exists $B:\T \to
\PSL(2,\R)$ analytic such that $B(x+\alpha)A(x)B(x)^{-1}$ is a triangular
matrix for every $x \in \T$?

\end{problem}
}

\subsubsection{Optimality} 

We do not claim optimality of the condition we impose on $\rho$.  Our proof
uses a single procedure, rational approximation, independent of the
Diophantine properties of $\alpha$.  This has the advantage of providing a
unified argument, but is clearly unoptimal in the
Diophantine case (it is actually surprising that the procedure works at
all in the case of bounded type).  More precise estimates should be
obtainable by interpolation arguments with classical KAM schemes.

\section{Notations and preliminaries}

\subsection{The fibered Lyapunov exponent.}\label{sec:2}  Given a cocycle $(\a,A)$, for $n \in \Z$, we denote the iterates of $(\a,A)$ by $(\a,A)^n=(n\a,A^{(n)}(\cdot))$ where for $n\geq 1$
$$\begin{cases}&A^{(n)}(\cdot)=A(\cdot+(n-1)\a)\cdots A(\cdot)\\
 &A^{(-n)}(\cdot)=A(\cdot-n\a)^{-1}\cdots A(\cdot-\a)^{-1}\end{cases}$$
We call fibered products of $(\a,A)$ the matrices $A^{(n)}(\cdot)$ for $n \in \N$.

The fibered Lyapunov exponent is defined as the limit
$$L(\a,A):= \lim_{|n| \to \infty} {1\over n} \int_{\theta \in {\T}} 
\log \|A^{(n)}(\theta)\| d\theta,$$
which by the subbadditive theorem always exists (similarly the limit when $n$ goes to $-\infty$ exists and is equal to $L(\a,A)$).

\subsection{The fibered rotation number.} \label{rot}

Assume that $A(\cdot):{\T}\to SL(2,{\R})$ is continuous and  homotopic to the identity; then the same is true for the map 
\begin{align*}F:{\T}\times{\mathbb S}^1&\to {\T}\times{\mathbb S}^1\\
(\theta,v)&\mapsto (\theta+\a,{A(\theta)v\over \|A(\theta)v\|}),\end{align*}
therefore $F$ admits  a continuous lift $ F:{\T}\times{\R}\to {\T}\times{\R}$ of the form $\tilde{F}(\theta,x)=(\theta+\a,x+f(\theta,x))$
 such that $f(\theta,x+1)=f(\theta,x)$ and $\pi (x+f(\theta,x))=A(\theta)\pi(x)/\|A(\theta)\pi(x)\|$, where $\pi: \R \to {\mathbb S}^1$, $\pi(x) = e^{i2 \pi x}:=(\cos(2\pi x),\sin(2\pi x))$. In order to simplify the terminology we shall say that $\tilde{F}$ is a lift for $(\a,A)$. The map $f$ is independent of the choice of  the lift up to the addition of a constant  integer $p\in{\Z}$. Following \cite{H} and \cite{JM}
we define the limit 
$$\lim_{n\to\pm\infty}{1\over n}\sum_{k=0}^{n-1} f({\tilde{F}}^k(\theta,x)),$$
that is independent of $(\theta,x)$ and where the convergence is uniform in $(\theta,x)$. The class of this number in 
$\T$, which is independent of the chosen lift, is called the {\it fibered rotation number} of $(\a,A)$ and denoted by $\rho(\a,A)$.
Moreover $\rho(\a,A)$ is continuous as a function of $A$ (with respect to the uniform topology on
$C^0(\T,SL(2,\R))$, naturally restricted to the subset of $A$ homotopic to
the identity).

\subsection{Continued fraction expansion.}\label{sec:2.3}

Define as usual for $0< \a <1$,
$$ a_0=0,\qquad \a_{0}=\a,$$
and inductively for $k\geq 1$,
$$a_k=[\a_{k-1}^{-1}],\qquad \a_k=\a_{k-1}^{-1}-a_k=G(\a_{k-1})=\{{1\over \a_{k-1}}\},$$
We define
\begin{align*}&p_0=0\qquad q_1=a_1\\
&q_0=1\qquad p_1=1,\end{align*}
and inductively,
\begin{equation} \label{3.1} \begin{cases}&p_k=a_kp_{k-1}+p_{k-2}\\&q_k=a_kq_{k-1}+q_{k-2}.\end{cases}\end{equation}


Recall that the sequence $(q_n)$  is  the sequence of best denominators  of $\a \in \R \setminus \Q$ since it satisifies 
\begin{equation*}  \forall 1 \leq k < q_n,\quad \norm{k\a} \geq \norm{q_{n-1}\a} \end{equation*}
and 
\begin{equation}
\norm{q_n \a } \leq {1 \over q_{n+1}} \label{red1}
\end{equation}
where we used the notation 
\begin{displaymath}
  \norm{x} = \|x\|_\T=\inf_{p \in \Z} | x - p|.
\end{displaymath}

\section{Growth of cocycles of rotations} \label{sec.den}

Let be given a number $\a \in \R \setminus \Q$. We will fix in the sequel a particular  subsequence ${(q_{n_k})}$ of the denominators of $\a$, that  we will denote  for simplicity by $(Q_k)$ and denote $(\bQ_k)$ the sequence $(q_{n_k+1})$.  
The goal of this section is to introduce a subsequence $(Q_k)$ for which we will have a nice control on the Birkhoff sums of real analytic functions above the irrational rotations of frequency $\a$. 
The properties required from our choice of the sequence $(Q_k)$ are summarized in the following statement, and are all what will be needed from this section in the sequel. The rest of Section \ref{sec.den} is devoted to the construction of the sequence $(Q_k)$ and can be skipped in a first reading.  

If $f\in C^0(\T,\R)$ we define its $k$-th Fourier coefficient by $\hat f(k)=\int_{\T}f(x)e^{-2\pi i kx}dx$ and  the $n$-th Birkhoff sum of $f$ over $x\mapsto x+\a$ by $S_{n}f=\sum_{k=0}^{n-1}f(\cdot+k\a)$.
We also introduce the notation $\eta_n=\eta/n^2$, for any $\eta>0$.

\begin{prop} \label{denjoy}
Given any $\eta \in(0,1)$, $h_*>0$ and $M>1$, there exists $C(h_*,\eta,M)>0$ such that for any irrational $\alpha$, there exists  a subsequence  $(Q_k)$ of denominators of $\a$ such that  $Q_0=1$, and for any $h\geq
h_*$ and any function  $\varphi \in
C^{\omega}_{h}(\T,\R)$, it holds for all $k>0$ and $h_k:=h(1-\eta_k)$
\begin{itemize} 
\item $Q_{k+1}\leq {\bQ}_k^{16M^{4}}$

\item   $\displaystyle{{\|S_{Q_k} \varphi-Q_k \widehat{\varphi}(0)\|}_{h_{k}} \leq  C {\|\varphi-\widehat{\varphi}(0)\|}_h ( {Q_k^{-M}} + {\bQ}_k^{-1+\frac{1}{M}} )}$.
\item for any $l \leq Q_{k+1}$, $\displaystyle{{\|S_{l} \varphi-l \widehat{\varphi}(0)\|}_{h_{k}} \leq C {\|\varphi-\widehat{\varphi}(0)\|}_h \left({\bQ_k}{Q_k^{-M}} +
{\bQ}_k^{\frac{1}{M}}\right)}$.
\end{itemize}
\end{prop}

\begin{defn} Let $0<\tA\leq \tB\leq \tC$. We say that the pair of denominators $(q_l,q_n)$ forms a
$\CD(\tA,\tB,\tC)$ bridge if 
\begin{itemize}
\item $q_{i+1}\leq q_i^{\tA}, \quad \forall i=l,\ldots,n-1$
\item $q_l^{\tC}\geq q_n\geq q_l^{\tB}$
\end{itemize} 
\end{defn} 

\begin{lem} \label{dioph.bridge} For any ${\tA}$ there exists  a subsequence $Q_k$ such that $Q_0=1$ and for each $k\geq 0$, 
$Q_{k+1}\leq \bQ_k^{{\tA}^4}$, and either $\bQ_k\geq Q_k^{\tA}$, or the pairs $(\bQ_{k-1},Q_{k})$ and $(Q_k,Q_{k+1})$ are both
$\CD({\tA},{\tA},{\tA}^3)$ bridges.
\end{lem}

\begin{pf} Assume the sequence $Q_l$ is constructed up to $k$. Let $q_n>Q_k$ (if it exists) be the smallest denominator such that
$q_{n+1}>q_n^{\tA}$. If $q_n\leq \bQ_k^{{\tA}^4}$ then we let $Q_{k+1}=q_n$. If $q_n\geq
\bQ_k^{{\tA}^4}$ (or if $q_n$ does not exist) it is possible to find $q_{n_0}:=\bQ_k$, $q_{n_1},q_{n_2},\ldots,q_{n_j}$ (or an infinite sequence $q_{n_1},q_{n_2},\ldots$) such that : $j\geq 2$, $q_{n_j}=q_n,$ and for each $0\leq i\leq j-1$, the pairs $(q_{n_i},q_{n_{i+1}})$  and $(q_{n_i+1},q_{n_{i+1}})$  are
$\CD({\tA},{\tA},{\tA}^3)$ bridges. We then let $Q_{k+1}:=q_{n_1}$, $Q_{k+2}:=q_{n_2},\ldots,Q_{k+j}:=q_{n_j}$ and it is straightforward to check that $Q_{k+1},Q_{k+2},\ldots,Q_{k+j-1}$ (or $Q_{k+1},Q_{k+2},\ldots$) satisfy the second condition of the lemma while $Q_{k+j}=q_n$, in case $q_n$ exists, satisfies the first one. 
\end{pf}

\medskip 

\noindent {\bf Proof of proposition \ref{denjoy}.} 
\begin{lem} \label{denjoy.el} Let $h_*,\eta,U>0$. There exists $C(h_*,\eta,U)>0$ such that for any $\varphi \in
C^{\omega}_{h}(\T,\R)$ with $h\geq h_*$ and for any irrational $\a$ and any pair $q_{s_1}\leq q_{s_2}$ of denominators of $\a$, we have if $\delta\geq {\rm max}(1/\sqrt{q_{s_1}},\eta/(10s_2^{2}))$ 
 $${\|S_{q_{s_2}} {\varphi}-q_{s_2}\widehat{{\varphi}}(0)\|}_{{h}(1-\delta)} \leq C {\|{\varphi}-\widehat{{\varphi}}(0)\|}_{{h}} \left( s_2^{4}\frac{q_{s_1}}{q_{s_2+1}} + \frac{q_{s_2}}{q_{s_1}^U}\right)$$
\end{lem}
\begin{pf}
Write $\varphi(z)=\sum_{l \in \Z} \widehat{\varphi}(l) e^{i2\pi lz}$.
 For any $l \in \Z^*$, we have that  \newline $|\widehat{\varphi}(l)| e^{2\pi |l|h} \leq  {\|\varphi-\widehat{\varphi}(0)\|}_h$. Now
$$S_{q_{s_2}}\varphi(z)-q_{s_2} \widehat{\varphi}(0) = \sum_{l \in \Z^*} \frac{1- e^{i2\pi q_{s_2} l \a}}{1- e^{i2\pi l \a}}    \widehat{\varphi}(l)e^{i2 \pi lz}$$
hence 
\begin{eqnarray*} {\|S_{q_{s_2}}\varphi -q_{s_2} \widehat{\varphi}(0)\|}_{h(1-\delta)} &\leq&  {\|\varphi-\widehat{\varphi}(0)\|}_h  \sum_{l \in \Z^*}   \left| \frac{1- e^{i2\pi q_{s_2} l \a}}{1- e^{i2\pi  l \a}}\right|   e^{-2 \pi h |l| \delta} \\
 &\leq&  2\pi  {\|\varphi-\widehat{\varphi}(0)\|}_h \frac{q_{s_1}}{q_{s_2+1}} \sum_{  0<|l| < q_{s_1}}   |l|  e^{-2 \pi h  |l| \delta} \\ &+&   {\|\varphi-\widehat{\varphi}(0)\|}_h  \sum_{ |l| \geq q_{s_1}} q_{s_2} e^{-2 \pi h |l| \delta} \end{eqnarray*}
where we have used the facts that ${\|l\a\|}_{\T} \geq \frac{1}{2q_{s_1}}$ for $0<|l|<q_{s_1}$  and ${\|q_{s_2} \a\|}_{\T} \leq {1}/{q_{s_2+1}}$.
The result now follows from the condition $\delta\geq {\rm max}(1/\sqrt{q_{s_1}},\eta/(10s_2^{2}))$ 
\end{pf}

In all the sequel we will let ${\tA}:=2M$ and $U=16M^4$ in Lemmas \ref{dioph.bridge} and \ref{denjoy.el}. 

A consequence of Lemma \ref{denjoy.el} is the following.
\begin{cor} \label{cor3} Given $h_*,M>0$, there exists $T_0(h_*,M)$ such that, for any  $\varphi \in
C^{\omega}_{h}(\T,\R)$ with $h\geq h_*$, 
and for any irrational $\a$ and any denominator $q_n$ of $\a$ such that $q_{n}\geq T_0$, we have
 \begin{enumerate} 
\item  If $q_{n+1}\geq q_n^{\tA}$ then 
$${\|S_{q_n} \varphi \|}_{h(1-\eta_n)}\leq  {\|\varphi-\widehat{\varphi}(0)\|}_h \frac{1}{q_n^M}.$$
\item If there exists $l$ such that $q_n^{1/{\tA}^3}\leq q_l\leq q_n^{1/{\tA}}$, then 
 $${\|S_{q_n} \varphi \|}_{h(1-\eta_n)}\leq {\|\varphi-\widehat{\varphi}(0)\|}_h\left({q_{n+1}^{-1+1/M}}+{q_n^{-M}}\right).$$
\end{enumerate}
 As a consequence we get for $Q_k\geq T_0$
 $$\displaystyle{{\|S_{Q_k} \varphi-Q_k \widehat{\varphi}(0)\|}_{h(1-\eta_k)} \leq  {\|\varphi-\widehat{\varphi}(0)\|}_h ( {Q_k^{-M}} + {\bQ}_k^{-1+\frac{1}{M}} )}$$

\end{cor}
\begin{pf} To prove $1$ apply Lemma \ref{denjoy.el} to $s_1=s_2=n$: indeed if $q_n$ is sufficiently large we have $\eta_n \geq 1/\sqrt{q_n}$, whence $${\|S_{q_{n}} {\varphi}-q_{n}\widehat{{\varphi}}(0)\|}_{h(1-\eta_n)} \leq C {\|\varphi-\widehat{\varphi}(0)\|}_h  \left( n^{4}\frac{q_{n}}{q_{n+1}} + \frac{1}{q_{n}^{U-1}}\right)$$ and $1$ follows if $q_n$ is sufficiently large form the fact that $q_{n+1}\geq
q_n^{\tA}$ and from the choices ${\tA}=2M$ and $U=16M^4$. 

To prove $2$ we let $s_1=l$, $s_2=n$, whence  for $q_n$ sufficiently large 
$${\|S_{q_{n}} {\varphi}-q_{n}\widehat{{\varphi}}(0)\|}_{h(1-\eta_n)} \leq C {\|\varphi-\widehat{\varphi}(0)\|}_h  \left( n^{4}\frac{q_{l}}{q_{n+1}} + \frac{q_n}{q_{l}^{U}}\right).$$ Since $q_n\leq
q_l^{{\tA}^3}$ implies that $n=\cO(\ln q_l)$ (the denominators grow at least geometrically) we get $2$ from
$q_n^{1/{\tA}^3}\leq q_l\leq q_n^{1/{\tA}}$ and the choices of ${\tA}$ and $U$.

The conclusion of the corollary follows from $1$ if $\bQ_k\geq Q_k^{\tA}$ and from $2$ if not since in this case
$Q_k^{1/{\tA}^3} \leq \bQ_{k-1} \leq Q_k^{1/{\tA}}$ as $(\bQ_{k-1},Q_k)$ is a
$\CD({\tA},{\tA},{\tA}^3)$ bridge (we also use the fact that $n_k\geq k$ in $Q_k=q_{n_k}$). 
\end{pf}

Another consequence of lemma \ref{denjoy.el} is the following.
\begin{cor} \label{cor4} Given $h_*,\eta,M>0$, there exists $T_0(h_*,\eta,M)$ such that, for any  $\varphi \in
C^{\omega}_{h}(\T,\R)$ with $h\geq h_*$, 
 and for any irrational $\a$ and
 for $Q_k\geq T_0$, we have for $m \leq Q_{k+1}$ 
\begin{equation} \label{sm} \displaystyle{{\|S_{m} \varphi- l \widehat{\varphi}(0)\|}_{h(1-\eta_k)} \leq {\|\varphi-\widehat{\varphi}(0)\|}_h \left( {\bQ_k}{Q_k^{-M}} + {\bQ}_k^{\frac{1}{M}}\right)}. \end{equation}
 
 \end{cor}

\begin{pf} We distinguish two cases.   

\noindent {\bf Case 1: $\bQ_k\geq Q_k^{\tA}$}. We let $u,v$ be such that $q_u=Q_k$ and $q_{v+1}=Q_{k+1}$ and write $m< Q_{k+1}$ as $m= \sum_{s=u}^{v} a_sq_s+b, a_s< q_{s+1}/q_s$, $b\leq q_u$ and apply lemma \ref{denjoy.el} to each $s$ with $s_1=s_2=s$ and obtain 
\begin{align*}  {\|S_{m} {\varphi}-m\widehat{{\varphi}}(0)\|}_{h(1-\eta_u)} &\leq C {\|\varphi-\widehat{\varphi}(0)\|}_h  \left( \sum_{s=u}^{v} \left( s^{4} + \frac{q_{s+1}}{q_{s}^{U}} \right) +q_u \right) \\
&\leq  C {\|\varphi-\widehat{\varphi}(0)\|}_h (\frac{\bQ_k}{Q_k^U}+q_u+v^5) \end{align*}
 where we used that for every $s\in [u+1,v]$ we have that $q_{s+1}\leq
q_s^{\tA}$ (by construction of the sequence $Q_k$). Since $\bQ_k\geq Q_k^{\tA}$ and $v=\cO(\ln \bQ_k)$ (because $Q_{k+1}\leq
\bQ_k^{{\tA}^4}$) we obtain (\ref{sm}).

\noindent {\bf Case 2: $\bQ_k \leq Q_k^{\tA}$}. By definition of the sequence $Q_k$ we then have that  $(\bQ_{k-1},Q_{k})$ and $(Q_k,Q_{k+1})$ are
$\CD({\tA},{\tA},{\tA}^3)$ bridges. We let $u=n_{k-1}+1$, $v=n_{k+1}-1$ and write
$m< Q_{k+1}$ as $m= \sum_{s=u}^{v} a_sq_s+b, a_s< q_{s+1}/q_s$, $b\leq q_u$ and apply lemma \ref{denjoy.el} to each $s$ with $s_1=s_2=s$ and obtain as before
\begin{align*}  {\|S_{m} {\varphi}-m\widehat{{\varphi}}(0)\|}_{h(1-\eta_u)} &\leq C {\|\varphi-\widehat{\varphi}(0)\|}_h  \left( \sum_{s=u}^{v} \left( s^{4} + \frac{q_{s+1}}{q_{s}^{U}} \right) +q_u \right) \\
&\leq  C {\|\varphi-\widehat{\varphi}(0)\|}_h (q_u+v^5) \end{align*}
because we have that for each $s\in[u,v]$ that $q_{s+1}\leq q_s^{\tA}$. From there we obtain (\ref{sm}) since  $v=\cO(\ln q_u)$  (from $q_v\leq
q_u^{{\tA}^4}$) and $q_u \leq Q_k^{1/{\tA}}$. 
\end{pf}
 



Proposition \ref{denjoy} is enclosed in the definition of the sequence $Q_k$ given in lemma \ref{dioph.bridge} and in the conclusions of corollary \ref{cor3} and corollary \ref{cor4}. $\hfill \Box$

\section{The inductive step.}

We denote by $\Omega(h)$ the set of cocycles $(\a,A)$ with $A\in C^{\omega}_{h}$, and $\Omega(h,\epss, \tau,\nu)\subset\Omega(h)$ the set of cocycles with $\a$ irrational and the fibered rotation number $\rho$ satisfying  
\begin{equation} \norm{q_n \rho} \geq {\epss} \max (q_n^{-\tau}, q_{n+1}^{-\nu}). \label{eq.dioph} \end{equation}
\subsection{} In this section we show that a close to constant cocycle, but not too close to $\pm
\id$, can be reduced to become essentially $e^{-\rho^2/\alpha}$-close to rotations, where $\rho$ is its rotation frequency in the base and $\rho$ its fibered rotation number.


For $\phi\in\C$ we denote by $R_{\phi}$ the matrix $\begin{pmatrix}\cos\phi& -\sin\phi\\ \sin\phi&\cos\phi\end{pmatrix}$.

\begin{prop}  \label{CT} 

For every $D\geq0,h_*>0$, there exist $\epsilon_0=\epsilon_0(h_*,D)>0$ and
$C_0=C_0(h_*,D)>0$ with the following properties.
Let $h>h_*$, $0<\delta<1$,
$\bar \alpha \in \R$, $\bar A \in C^\omega_h(\T,\SL(2,\R))$,
$\bar \phi \in C^\omega_h(\T,\R)$ satisfy
$\| \bar \phi-\hat{\bar{\phi}}(0) \|_h \leq D$,
\be \label{eq.rho}
\rho^{-1}=\|(R_{2 \bar \phi}-\id)^{-1}\|_h<\epsilon_0^{-1/4},
\ee
\be
\|R_{-\bar \phi} \bar A-\id\|_h<\epsilon_0.
\ee
Then:
\begin{enumerate}
\item There exist
$B \in C^\omega_{e^{-\delta/3} h}(\T,\SL(2,\R))$ with
\begin{equation} \label {B1}
\|B-\id\|_{e^{-\delta/3} h} \leq \frac {C_0} {\rho^2}   \|R_{-\bar \phi} \bar A-\id\|_h 
\end{equation}
and $\tilde {\bar
\phi} \in C^\omega_{e^{-\delta/3} h}(\T,\R)$ such that, letting $\tilde
{\bar A}(x)=B(x+\bar \alpha) A(x) B(x)^{-1}$, we have
\begin{equation} \label {xi1}
\|R_{-\tilde {\bar \phi}} \tilde {\bar A}-\id\|_{e^{-\delta/3} h}
\leq C_0 e^{-\frac {h \delta \rho^2} {C_0 |\bar \alpha|}}
\|R_{-\bar \phi} \bar A-\id\|_h.
\end{equation}
\item If $\alpha \in \R$ and $A \in C^\omega_h(\T,\SL(2,\R))$ satisfy
$\|A\|_h \leq D$ and
$A(x+\bar \alpha) \bar A(x)=\bar A(x+\alpha) A(x)$ (i.e.,
$(\alpha,A)$ commutes with $(\bar \alpha,\bar A)$) then there
exists $\tilde \phi \in C^\omega_{e^{-\delta} h}(\T,\R)$ such that,
letting $\tilde A(x)=B(x+\alpha)A(x)B(x)^{-1}$ we have
\begin{equation} \label {xx1}
\|R_{-\tilde \phi} \tilde A-\id\|_{e^{-\delta} h} \leq
C_0 e^{-\frac {h \delta \rho^2} {C_0 |\bar \alpha|}}.
\end{equation}
\end{enumerate}
\end{prop}

\begin{pf}

We will need a preliminary lemma.

\begin{lemma} \label{elliptic}

For every $D>0$, there exists $C>0$,
$\epsilon>0$ with the following properties.
Let $W \subset \SL(2,\C) \times \C$ be the set of all
$(A,\theta)$ such that $\|A\|<2 D$ and
$\|R_\theta^{-1} A-\id\|<\epsilon \max \{1,\|R_{2 \theta}-\id\|^2\}$.
There exists a real symmetric
holomorphic function $F:W \to \SL(2,\C)$
such that $B=F(A,\theta)$ satisfies
$BAB^{-1}=R_{\theta'}$ and $\|B-\id\|<C \frac {\|R_\theta^{-1} A-\id\|}
{\|R_{2 \theta}-\id\|^2}$.

\end{lemma}

\begin{pf}

Let $A^{(0)}=A$, $\theta^{(0)}=\theta$.  Assuming $A^{(n)}$ and
$\theta^{(n)}$ defined and satisfying the inductive estimates
$|\theta^{(n)}-\theta| \leq (1-2^{-n})
\epsilon_0^{1/2} \|2 \theta\|^2_{\C/\Z}$ and
$\|A^{(n)}-R_{\theta^{(n)}}\|<\epsilon_0^{n/2} \|A-R_\theta\|$, define
$v^{(n)} \in \sl(2,\C)$ small such that
$A^{(n)}=e^{v^{(n)}} R_{\theta^{(n)}}$.  Let $v^{(n)}=\left ( \bm
x^{(n)} & y^{(n)}-2 \pi z^{(n)} \\ y^{(n)}+2 \pi z^{(n)} & -x^{(n)}
\em \right )$.  Let $\left (\bm \tilde x^{(n)} \\ \tilde y^{(n)} \em \right
)=(R_{2 \theta^{(n)}}-\id)^{-1} \left ( \bm x^{(n)} \\ y^{(n)} \em \right )$
and let $w^{(n)}=\left ( \bm \tilde x^{(n)} & \tilde y^{(n)} \\
\tilde y^{(n)} & -\tilde x^{(n)} \em \right )$.  Let $A^{(n+1)}=e^{w^{(n)}}
A^{(n)} e^{-w^{(n)}}$ and $\theta^{(n+1)}=\theta^{(n)}+z^{(n)}$.
Then $A^{(n+1)}$ and $\theta^{(n+1)}$ satisfy the inductive estimates as
well, so all the objects can be defined for every $n$.
Then one can take $B=\lim_{n \to \infty} e^{w^{(n-1)}} \cdots e^{w^{(0)}}$.
\end{pf}

In what follows we fix $C_2>10$ and $C(D)$ as in Lemma \ref{elliptic} and set $N=[\frac {\delta h \rho^2} {C_1 |\bar \alpha|}]$, where $C_1 \gg CC_2$ is a
sufficiently large constant.  Let $h_j=e^{-\delta \frac {j} {3 N}} h$,
$j \geq 0$.

\begin{claim}

If $\epsilon_0$ is sufficiently small, there exist sequences
$B_i,A_i \in C^\omega_{h_i}(\T,\SL(2,\R))$,
$\xi_i \in C^\omega_{h_i}(\T,\M(2,\R))$ and
$\phi_i \in C^\omega_{h_i}(\T,\R)$, $0 \leq i \leq N$,
such that
\begin{enumerate}
\item $A_0=\bar A$, $\phi_0=\bar \phi$,
\item $R_{\phi_i}=B_{i-1} A_{i-1} B_{i-1}^{-1}$ for $1 \leq i \leq N$
\item $A_i(x)=B_{i-1}(x+\bar\alpha) B_{i-1}(x)^{-1} R_{\phi_i(x)}$
for $1 \leq i \leq N$,
\item $\xi_i=R_{\phi_i}^{-1} A_i-\id$ for $0 \leq i \leq N$,
\item $\|B_i-\id\|_{h_i} \leq \frac {C} {\rho^2} \|\xi_i\|_{h_i}$ for
$0 \leq i \leq N-1$,
\item $\|R_{\phi_i}\|_{h_i} \leq (1+\frac {i} {N}) D$ for $0 \leq i \leq N$,
\item $\|(R_{2 \phi_i}-\id)^{-1}\|_{h_i} \leq (1+\frac {i} {N}) \rho^{-1}$ for
$0 \leq i \leq N$.
\item $\|\xi_i\|_{h_i} \leq \frac {1} {C_2}
\|\xi_{i-1}\|_{h_{i-1}}$ for $1 \leq i \leq N$.
\end{enumerate}

\end{claim}

\begin{pf}

We proceed by induction on $i$.  Assume that $B_0,...,B_{i-1}$ have been
already defined and hence $A_0,...,A_i$, $\phi_0,...,\phi_i$ and
$\xi_0,...,\xi_i$ are automatically defined by properties (1-4).  Assume
that estimates (5-8) have been also established for the objects so far
defined.

Let $F$ be as in the previous lemma.  We want to define
$B_i(z)=F(A_i(z),\phi_i(z))$.  For this, we must check that
$(A_i(z),\phi_i(z))$ is in the domain of $F$, which is in fact a consequence
of estimates (6-8).  Thus $B_i$ is well defined, and satisfies estimate (5)
as well by the previous lemma.  If $i<N$, we now define $A_{i+1}$,
$\phi_{i+1}$ and $\xi_{i+1}$ by (2-4).  Estimates (6) and (7) immediately
follow from (2) and (5) (from (8) we have that $\|\xi_{i}\|_{h_{i}} \leq \eps_0$ and (\ref{eq.rho}) then yields  $\|\xi_{i}\|_{h_{i}}/\rho^2 \leq \eps_0^{1/2}$).

We obviously have
$\|\xi_{i+1}(x)\| \leq \|B_i(x)\|^3 \|A_i(x)\|^2 \|B_i(x+\bar
\alpha)-B_i(x)\|$.  Thus, estimating the derivative of $B_i-\id$ with the
Cauchy formula, we get $\|\xi_{i+1}\|_{h_{i+1}} \leq \frac {C N \bar\a} {\delta h}
\|B_i\|_{h_{i+1}}^3 \|A_i\|_{h_{i+1}}^2 \|B_i-\id\|_{h_i}$, which immediately
implies (8) since $N=[\frac {\delta h \rho^2} {C_1 |\bar \alpha|}]$ with $C_1 \gg CC_2$. 
\end{pf}

We now set $B=B_{N-1} \cdots B_0$, $\tilde {\bar \phi}=\phi_N$, $\tilde {\bar A}=A_{N}$
so that (1-4) give $\tilde {\bar A}(x)=B(x+\bar \alpha)A(x)B(x)^{-1}=
R_{\tilde {\bar \phi}}(\id+\xi_N)$.
By (5) and (8) we have
(\ref {B1}), while (8) gives (\ref {xi1}).
This proves the first statement.

For convenience of notation, let us write $\phi=\phi_N$ and $\xi=\xi_N$ (so
that $\tilde {\bar A}=R_\phi (\id+\xi)$).

\def\cP{{\mathcal {P}}}
\def\cQ{{\mathcal {Q}}}

\begin{claim}

If $C_3(D)$ is sufficiently large and  if $\epsilon_0>0$ is sufficiently small and
$|\bar \alpha|<C_3^{-1} \delta h \rho^2$ then
\be
\|(R_{\phi(x+\alpha)+\phi(x)}-\id)^{-1}\|_{e^{-2
\delta/3} h}<3 \rho^{-1}.
\ee

\end{claim}

\begin{pf}

By the Cauchy formula, $\|\tilde A(x+\bar
\alpha)-\tilde A(x)\|_{e^{-2 \delta/3} h} \leq
\frac {C} {\delta h} |\bar \alpha| \ll \rho^2$.
The commutation relation gives $\|\tilde A(x) R_{\phi(x)} \tilde
A(x)^{-1}-R_{\phi(x+\alpha)}\|_{e^{-2 \delta/3} h} \leq \frac {C} {\delta h}
|\bar \alpha|+C \epsilon_0$, hence
$\|\cos 2 \pi \phi(x+\alpha)-\cos 2 \pi \phi(x)\|_{e^{-2 \delta/3} h} \ll \rho^2$.  Thus for each $x$ with $|\Im x|<e^{-2 \delta/3}
h$, one of $\|R_{\phi({x+\alpha}) \pm \phi(x)}-\id\| \ll \rho$.  Since $\|(R_{2 \phi(x)}-\id)^{-1}\| \leq
\rho^{-1}$, one of $\|R_{\phi({x+\alpha}) \pm \phi(x)}-\id\|_{e^{-2 \delta/3}
h} \ll \rho$.
Clearly $\|R_{\phi(x+\alpha)+\phi(x)}-\id\|_{e^{-2 \delta/3} h}$ is of order
at least $\rho$, since
\be
\inf_{x \in \T} \|R_{\phi(x+\alpha)+\phi(x)}-\id\| \geq
\inf_{x \in \T} \|R_{2 \phi(x)}-\id\|,
\ee
and we must have
$\|R_{\phi(x+\alpha)-\phi(x)}-\id\|_{e^{-2 \delta/3} h} \ll \rho$.  Since $\|(R_{2 \phi}-\id)^{-1}\|_{e^{-\delta/3} h}$
is at most $2 \rho^{-1}$, the claim follows.
\end{pf}

Notice that the second statement holds trivially unless $|\bar \alpha| \ll \delta h \rho^2$, so we will assume from now on that
the estimate of the claim holds.

Let $J=\left (\bm 0 & -1 \\ 1 & 0 \em \right )$.  For $M \in \M(2,\C)$,
let
$\cQ(M)=\frac {M+JMJ} {2}$.  Then $M-\cQ(M)$ is of the form $\left (\bm a &
b \\ -b & a \em \right )$, so that it commutes with all $R_\theta$, $\theta
\in \C$.
Notice that for $\theta \in \C$ we have
\be
R_{-\theta} \cQ(M)=\cQ(M R_\theta), \quad
R_\theta \cQ(M)=\cQ(R_\theta M),
\ee
so for $\theta_1,\theta_2 \in \C$ we have
\be \label {Q}
(R_{-\theta_1}-R_{\theta_2}) \cQ(M)=\cQ(M R_{\theta_1}-R_{\theta_2} M),
\ee

Let $L(x)=\cQ(\tilde A(x))$, $L_1(x)=L(x+\bar \alpha)-L(x)$ and
$L_2(x)=\cQ(\tilde A(x) R_{\phi(x)}-R_{\phi(x+\alpha)} \tilde A(x))$.
Then for $2 N \leq n \leq 3N-1$ the Cauchy formula gives
\be
\|L_1\|_{h_{n+1}} \leq \frac {C |\bar \alpha|} {h_n-h_{n+1}}
\|L\|_{h_n},
\ee
for some absolute $C$, a direct computation bounds $\|L_2(x)\|$ by
\be
\|L_1(x)\| \|R_{\phi(x)}\|+\|\tilde
A(x+\bar \alpha)\| \|\tilde {\bar A}(x)-R_{\phi(x)}\|+
\|\tilde {\bar A}(x+\alpha)-R_{\phi(x+\alpha)}\| \|\tilde A(x)\|
\ee
(here we have used the commutation relation), giving
\be
\|L_2\|_{h_{n+1}} \leq C (\|L_1\|_{h_{n+1}}+\|\xi\|_{h_{n+1}}),
\ee
where $C$ only depends on $D$, and we have used that $\epsilon_0$ is small,
and (\ref {Q}) gives
\be
\|L\|_{h_{n+1}} \leq \|(R_{-\phi(x)}-R_{\phi(x+\alpha)})^{-1}\|_{h_{n+1}}
\|L_2\|_{h_{n+1}}.
\ee
Thus
\be
\|L\|_{h_{n+1}} \leq C \frac {\rho^{-1} |\bar \alpha|}
{h_n-h_{n+1}} \|L\|_{h_n}+C \rho^{-1} \|\xi\|_{h_N},
\ee
so that
\be \label {l}
\|L\|_{h_{3N}} \leq C e^{-\frac {h \delta} {C \rho |\bar \alpha|}}+C
\rho^{-1} \|\xi\|_{h_N} \leq
C e^{-\frac {h \delta \rho^2} {C |\bar \alpha|}}
\ee
where the second inequality uses that $\|\xi\|_{h_N} \leq
C \rho^4 e^{-\frac {h \delta \rho^2} {C |\bar \alpha|}}$.

Since we are assuming that $\bar \alpha \ll \delta h
\rho^2$, $L=\cQ(\tilde A)$ is
small so that $\frac {\tilde A-L} {(\det (\tilde A-L))^{1/2}}$ defines a
map in $C^\omega_{e^{-\delta} h}(\T,\SO(2,\R))$ which is $C
\|L\|_{e^{-\delta} h}$ close to $\tilde A$.

Notice that $\tilde A$ is itself homotopic to a constant since
$(\alpha,\tilde A)$ and $(\bar \alpha,\tilde {\bar A})$ commute
and $\tilde {\bar A}$ is homotopic to $R_\phi$ which is homotopic to a
constant (since $\phi$ takes values on $\R$).  Thus we can write $\frac
{\tilde A-L} {(\det (\tilde A-L))^{1/2}}=R_{\tilde \phi}$ for some function
$\tilde \phi \in C^\omega_{e^{-\delta} h}(\T,\R)$.  We then get
$\|R_{-\tilde \phi} \tilde A-\id\|_{e^{-\delta} h} \leq C \|L\|_{e^{-\delta}
h}$, which together with (\ref {l}) gives the second statement.
\end{pf}

\subsection{}
The following lemma shows that if a cocycle in $\Omega(h,\epss, \tau,\nu)$ is sufficiently close to rotations, then it is possible to iterate it and use the Diophantine property on the fibered rotation number and the bounds obtained for the growth of  Birkhoff sums involving the special denominators $Q_k$  to end up with  the required conditions for the reduction step of Proposition \ref{CT}.

In the following statement and all through  the rest of the paper we assume $\epss, \tau,\nu$ given and let 
 \begin{equation} M= {\rm max}(4\btau,\frac{2}{1-2\bnu}), \quad \btau=\tau+1, \quad \bnu = \frac{1}{2}(\nu+\frac{1}{2}) \label{defM} \end{equation}
and define the sequence $(Q_k)$ as in proposition \ref{denjoy} and set 
$$U_k= e^{-\bQ_kQ_k^{-b}-\bQ_k^{a}}, \quad a=\frac{2}{M} \quad b=\frac{M}{2}.$$
Recall that for $\eta>0$, $\eta_k=\eta/k^2$.

\begin{lem} \label{Aq}  Fix $D,\eta,h_*>0$. There exists $J(h_*,D,\eta,\epss,\nu,\tau)$ such that if $\bQ_k\geq J$ and if $(\a,A) \in \Omega(h,\epss, \tau,\nu)$  for some $h> h_*$  is such that $A = R_{\varphi}(\id +\xi)$  with ${\| \varphi-\widehat{\varphi}(0) \|}_h \leq D$ and ${\|\xi \|}_h \leq U_k$,
then we can express $A^{(Q_{k+1})}$ as $R_{\varphi^{(Q_{k+1})}}(\id+\xi^{(Q_{k+1})})$ with $\varphi^{(Q_{k+1})}:=S_{Q_{k+1}} \varphi  $ and

\begin{enumerate}
\item $ \| {\varphi^{(Q_{k+1})}} -\widehat{\varphi^{(Q_{k+1})}}(0)
\|_{h(1-\eta_{k+1})} \leq D$
\item  $\|(R_{2 \varphi^{(Q_{k+1})}}-\id)^{-1}\|_{h(1-\eta_{k+1})}<\rho_k, \ \ \rho_k= \frac{\varepsilon}{4(Q_{k+1}^{-\tau}+\bQ_{k+1}^{-\nu})}<U_k^{-\frac{1}{10}}$
\item ${\| \xi^{(Q_{k+1})} \|}_{h(1-\eta_{k+1})} \leq U_k^{\frac{1}{2}}$
\end{enumerate}

\end{lem}
\begin{pf} Given $C=C(h_*,\eta,M)$ of Proposition \ref{denjoy}, let $J>0$ be such that if $\bQ_k \geq J$ then 
\begin{align}
\label{J1} C D ({\bQ_k}{Q_k^{-M}} + {\bQ}_k^{\frac{1}{M}})&\leq\frac{1}{100} ( \bQ_k Q_k^{-b} + {\bQ}_k^{a}) \\
\label{J2} C  D ( Q_{k+1}^{-M} + {\bQ}_{k+1}^{-1+\frac{1}{M}})&\leq\frac{\epss}{100}(Q_{k+1}^{-\tau}+\bQ_{k+1}^{-\nu}) \\
\label{J3} U_k^{\frac{1}{10}}&\leq\frac{\epss}{100}Q_{k+1}^{-\tau} \leq \frac{1}{100}   \norm{Q_{k+1}\rho}   \end{align}
the last requirement being possible since by definition of the sequence $(Q_k)$ we have that $Q_{k+1}\leq \bQ_k^{16M^4}$. 
It is henceforth assumed that $\bQ_k\geq J$.
From Proposition \ref{denjoy} and (\ref{J1}), it follows  that for all $l\leq Q_{k+1}$, 
\begin{equation} {\|S_l \varphi-l \widehat{\varphi}(0)\|}_{h(1-\eta_{k+1})} \leq  \frac{1}{100} ( \bQ_k Q_k^{-b} + {\bQ}_k^{a}).  \label{eq.rot} \end{equation}


To prove (3)  we need the following 
straightforward claim on the growth of matrix products that can be proved following the lines of  Lemma 3.1. of \cite{AK1} 
\medskip

\begin{claim}

 We have that 
$$M_l (\id+\xi_l)\ldots M_0 (\id+\xi_0)=M^{(l)}(\id+\xi^{(l)})$$
with $M^{(l)}=M_l\ldots M_0$ and $\xi^{(l)}$ satisfying
$$\|\xi^{(l)}\|\leq e^{\sum_{k=0}^{l} \|M^{(k)}\|^2 \xi_{k}}-1$$

\end{claim}

To prove (3) apply the claim to $M_k=R_{\varphi(x+k\a)}$ and $\xi_k=\xi(x+k\a)$ and observe that $\|R_\theta\|\leq e^{|{\rm Im} \theta |}$, and use (\ref{eq.rot}).

On the other hand, Proposition \ref{denjoy} together with (\ref{J2}) and the Diophantine condition (\ref{eq.dioph}) on the fibered rotation number  imply that 
 \begin{equation} \label{422} {\| \varphi^{(Q_{k+1})}-Q_{k+1}
\widehat{\varphi}(0)\|}_{h(1-\eta_{k+1})}  \leq   \frac{\epss}{100}(Q_{k+1}^{-\tau}+\bQ_{k+1}^{-\nu}) \leq  \frac{1}{50} \norm{Q_{k+1}\rho}.\end{equation}
 Clearly (\ref{422})  implies (1).
 
 It comes from (\ref{422}) and  (\ref{J3}) and (3) that 
 $${\| \varphi^{(Q_{k+1})}-Q_{k+1} \widehat{\varphi}(0)\|}_{h(1-\eta_{k+1})}+{\| \xi^{(Q_{k+1})}
\|}_{h(1-\eta_{k+1})} \leq \frac{1}{25} \norm{Q_{k+1}\rho}.$$
But the fibered rotation number of the cocycle $(Q_{k+1}\a,A^{(Q_{k+1})})$ is $Q_{k+1}\rho$, hence 
$${\| \varphi^{(Q_{k+1})}-Q_{k+1} \rho\|}_{h(1-\eta_{k+1})} \leq  \frac{1}{10} \norm{Q_{k+1}\rho}$$
so that using  (\ref{J3}) again and the Diophantine property on the fibered rotation number $\rho$ we get (2) and the proof of the lemma is over.

\end{pf}

\subsection{}
As a consequence of Proposition \ref{CT} and  Lemma \ref{Aq} we get the following inductive step, fundamental in our reduction. 
\begin{prop} \label{prop.inductive}  Given $h_*,\eta,D>0$, there exists $T(h_*,D,\eta,\epss,\nu,\tau)$ such that if $k\geq 1$, $\bQ_k\geq T$ and if
$(\alpha,A_k) \in  \Omega(h,\epss,\tau,\nu)$ for some  $h\geq h_*$, can be
written as $A_k = R_{\varphi_k}(\id +\xi_k)$ with
 \begin{itemize}
 \item ${\|\varphi_k - \widehat{\varphi}_k(0)\|}_{h} \leq D-\eta_k$
 \item $\displaystyle{{\| \xi_k \|}_{h} \leq U_k}$ \end{itemize}
 then,  if we denote  $h_{k+1} :=h (1-\eta_{k+1})^2$, 
 there exists $B_k$ with  ${\|B_k-\id\|}_{h_{k+1}} \leq  U_k^{\frac{1}{4}}$  such that
$A_{k+1}(x):=B_k(x+\a)A_k(x)B_k(x)^{-1}$ can be written as
$A_{k+1} = R_{\varphi_{k+1}}(\id +\xi_{k+1})$ with
 \begin{itemize}
 \item ${\|\varphi_{k+1}  - \widehat{\varphi}_{k+1}(0)\|}_{h_{k+1}} \leq
 D-\eta_{k+1}$
  \item ${\| \xi_{k+1} \|}_{h_{k+1}} \leq  U_{k+1}$
 \end{itemize}
 \end{prop}
 
\begin{pf} Let $\epsilon_0(h_*/2,2D)$ be as in Proposition \ref{CT}. 
Assume $k$ is such that
 $\bQ_k \geq J$ where $J(h_*,D,\eta,\epss,\nu,\tau)$ is given by Lemma \ref{Aq} and  $U_k \leq \epsilon_0(h_*/2,1)^2$. 
 

By Lemma \ref{Aq}  we can apply Proposition \ref{CT} to the cocycle $(\bar \a,\bar A)$ with $\bar\a=Q_{k+1}\a$    and $\bar A=A_k^{(Q_{k+1})}$, and to the cocycle $(\a,A_k)$    that commutes with  $(\bar \a,\bar A)$. We  thus get  a conjugacy $B_k$ such that
${\|B_k-\id\|}_{h(1-\eta_{k+1})} \leq C_0 U_k^{1/2-1/5}$ (see (\ref{B1}))  
while (\ref{xx1}) yields that $A_{k+1}(x)=B_k(x+\a) A_k(x) B_k(x)^{-1}$ can be expressed as $A_{k+1} =  R_{\varphi_{k+1}}(\id +\xi_{k+1})$ with 
 $${\| {\xi_{k+1}}  \|}_{h_{k+1}} \leq C_0 e^{-\frac{h(1-\eta_{k+1})
\eta_{k+1} \rho_k^2}{ C_0 \norm{Q_{k+1}\a}}}$$
with  $\rho_k= \frac{\varepsilon}{4(Q_{k+1}^{-\tau}+\bQ_{k+1}^{-\nu})}$.

Since  $\norm{Q_{k+1}\a}\leq 1/\bQ_{k+1}$, (\ref{defM}) implies that $ {\bQ_{k+1}^{2/M}+\bQ_{k+1}Q_{k+1}^{-M/2}} = {o} (
\eta_{k+1} \rho_{k}^2/\norm{Q_{k+1}\a} )$ as $\bQ_k$ tends to infinity. Hence, there exists $T(h_*,D,\eta,\epss,\nu,\tau)$, such that for
 $\bQ_k\geq T$ we have    ${\| \xi_{k+1} \|}_{h_{k+1}} \leq  U_{k+1}$, while the bound on $\varphi_{k+1}$ follows immediately from  the bounds  on $\varphi_k$ and $B_k$.
\end{pf}

\begin{cor} \label{corfin} Let $\epss,\tau>0$, $h>2 \epss$, $0<\nu<1/2$, and let
$T(h-\epss,\epss,\epss/10,\epss,\nu,\tau)$ be as in Proposition \ref{prop.inductive}.
Suppose that  $(\a,A') \in \Omega(h-\epss,\epss,\tau,\nu)$ with $A' =R_{\varphi'}(\id +\xi')$ satisfying
\be
{\|\varphi' - \widehat{\varphi}'(0)\|}_{h-\epss/2} \leq \frac{\epss}{2}
\ee
\be
\| \xi' \|_{h-\epss/2} \leq U_{n_0}
\ee
with $n_0$ such that $\bQ_{n_0} \geq T$. Then there exist $B:\Delta_{h-\epss} \to {\rm SL}(2,\C)$
and $\varphi:\Delta_{h-\epss} \to \C$ such that
\begin{enumerate}
\item $B$ and $\varphi$ are real-symmetric,
\item $\|B-\id\|_{h-\epss}<\epss/2$ and
$\|\varphi- \widehat{\varphi}(0)\|_{h-\epss}<\epss$,
\item $B(x+\alpha)A'(x)B(x)^{-1}=R_{\varphi(x)}$.
\end{enumerate}
\end{cor}

\section{Proof of Theorem \ref{theorem.cd}.}

Let $h,\epss,\tau>0$ and $0<\nu<1/2$ be given (we assume that $\varepsilon<h/2$).  

\comm{
The following is an immediate consequence of an inductive application of  Proposition \ref{prop.inductive}. 
\begin{cor} \label{corfin} Take $\eta=\epss/10$ and let $T(h-\epss,\epss,\eta,\epss,\nu,\tau)$ be as in Proposition \ref{prop.inductive}. Suppose that  $(\a,A') \in \Omega(h-\epss,\epss,\tau,\nu)$ with $A' =R_{\varphi'}(\id +\xi')$ satisfying
\begin{itemize}
\item[(H1)] ${\|\varphi' - \widehat{\varphi}'(0)\|}_{h-\epss/2} \leq \frac{\epss}{2}$
\item[(H2)] $\displaystyle{{\| \xi' \|}_{h-\epss/2} \leq U_{n_0}}$
\end{itemize}
with $n_0$ such that $\bQ_{n_0} \geq T$. Then there exists $B:\Delta_{h-\epss} \to {\rm SL}(2,\C)$
and $\varphi:\Delta_{h-\epss} \to \C$ such that
\begin{enumerate}
\item $B$ and $\varphi$ are real-symmetric,
\item $\|B-\id\|_{h-\epss}<\epss/2$ and
$\|\varphi- \widehat{\varphi}(0)\|_{h-\epss}<\epss$,
\item $B(x+\alpha)A'(x)B(x)^{-1}=R_{\varphi(x)}$.
\end{enumerate}
\end{cor}
}

In order to obtain Theorem \ref{theorem.cd}, it is sufficient to show that there exists $B'$ with
\be \label {B'}
{\| B' -\id\|}_{h-\epss/2} \leq \epss/10
\ee
and $(\alpha,A')$ satisfying the conditions of Corollary \ref{corfin}
such that $A'(x)=B'(x+\a)A(x)B'(x)^{-1}$.

We first claim that there exists $n_0$ such that $Q_{n_0} \leq T^{{\tA}^4}$ while $\bQ_{n_0}\geq T$. Indeed, 
let $m_0$ be such that $Q_{m_0} \leq T \leq Q_{m_0+1}$. There are two
possibilities: either $\bQ_{m_0}\geq T$ (and $n_0=m_0$ satisfies the claim)
or $\bQ_{m_0}\leq T$ and by definiton of the sequence $(Q_k)$ it then holds that
$Q_{m_0+1}\leq T^{{\tA}^4}$ while by definition of $m_0$ it holds that $\bQ_{m_0+1}\geq T$
(so $n_0=m_0+1$ satisfies the claim).

Let $\epsilon_0(h-\epss,\epss)$ be as in Proposition \ref{CT}. Define $\epsilon_1:={\rm min}(\epsilon_0,(\frac{\epss}{4}Q_{n_0}^{-\tau})^4)$.
Since $Q_{n_0} \leq T^{{\tA}^4}$, by taking $\epsilon$ sufficiently small
(depending on $h,\epss,\tau,\nu$) then the hypothesis that $A$ is
$\epsilon$-close to a constant rotation through $\Delta_h$ implies that
$\xi_0=R_{-Q_{n_0} \rho} A^{(Q_{n_0})}$ satisfies $\|\xi_0\|_h<\epsilon_1$,
(where $\rho$ is the fibered rotation number of the cocycle $(\a,A)$).

Now, if $(\a,A) \in \Omega(h,\epss, \tau,\nu)$, then $\norm{2 Q_{n_0} \rho} \geq  \epss Q_{n_0}^{-\tau}   \geq 4\eps_1^{1/4}$.
This  implies that  ${\|(R_{2 Q_{n_0} \rho} -\id)^{-1} \|}  \leq
2/\norm{2 Q_{n_0} \rho} \leq {\|\xi_{0}\|}_{h}^{-1/4}/2$. We are thus in position to  apply Proposition \ref{CT}  to
$(\overline \alpha,\overline A)=(\alpha,A)^{Q_{n_0}}$ and get a
conjugacy $B'$ such that ${\|B'-\id\|}_{h_{0}-\epss/4} \leq C_0 \epsilon_0^{\frac{1}{2}}$
and $A'(x):=B'(x+\alpha)A(x)B'(x)^{-1}$ can be written as $A'=R_{\varphi'}(\id+\xi')$
where  $${\|\xi'\|}_{h_{0}-\epss/2} \leq C_0  e^{-\frac{h \epss  {\norm{Q_{n_0} \rho}}^2}{16C_0\norm{Q_{n_0}\alpha}}}$$ 
and as shown in the proof of Proposition \ref{prop.inductive} it follows from the fact that $\bQ_{n_0} \geq T$ that 
${\| {\xi'}  \|}_{h-\epss/2} \leq U_{n_0}$.

Moreover, it follows from the bound on $B'$ that  ${\|\varphi' - \widehat{\varphi'}(0)\|}_{h-\epss/2} \leq \frac{\epss}{2}$. 
\qed


\end{document}